\documentclass[11pt]{article}
\usepackage{amssymb}
\usepackage{epsfig}
\usepackage{epstopdf}
\usepackage[usenames]{color}
\usepackage[numbers]{natbib}
\def\citet{\cite}
\usepackage{amssymb}

\newtheorem{theorem}{Theorem}[section]

\newtheorem{lemma}[theorem]{Lemma}
\newtheorem{corollary}[theorem]{Corollary}

\newtheorem{remark}[theorem]{Remark}
\newtheorem{example}[theorem]{Example}

\def\a{\alpha}

\def\O{\Omega}

\def\F{{\cal F}}

\def\esssup{\mathop{\rm ess\, sup}}

\def\R{{\bf R}}

\def\P{{\bf P}}

\def\s{\delta}
\def\g{\gamma}
\def\C{{\bf C}}

\def\ww{\widetilde}

\def\s{\sigma}

\def\GG{{\cal G}}

\newcommand{\be}{\begin{equation}}
\newcommand{\ee}{\end{equation}}
\newcommand{\bd}{\begin{displaymath}}
\newcommand{\ed}{\end{displaymath}}
\newcommand{\ba}{\begin{array}{ll}}
\newcommand{\ea}{\end{array}}
\newcommand{\baa}{\begin{eqnarray}}
\newcommand{\eaa}{\end{eqnarray}}
\newcommand{\baaa}{\begin{eqnarray*}}
\newcommand{\eaaa}{\end{eqnarray*}}
\font\sm=cmr10

\def\Re{{\rm Re\,}}
\def\Im{{\rm Im\,}}

\def\o{\omega}

\def\mm0{m_{\scriptscriptstyle 0}}
\def\m1{m_{\scriptscriptstyle 1}}

\title{On stochastic integrals with controlled growth  of their containing range}
\author{
Nikolai Dokuchaev\\ \ {\sm Department of Mathematics \& Statistics,
Curtin University,}\\
{\sm  GPO Box U1987, Perth, 6845 Western Australia}} 
\begin{document}
\maketitle
\begin{abstract} This short note suggests special examples of stochastic It\^o integrals with controlled growth  of their containing range.  The integrands for this integrals are presented
explicitly. The construction does not involve neither stopping times nor forecasting or calculation of the conditional expectations of a contingent claim.

{\bf Key words}: stochastic integrals, It\^o calculus, containing range

{\bf Mathematical Subject Classification (2010):} 
65C30, 
		65C50, 
		65C60 
\end{abstract}
\section{Introduction}
The paper considers stochastic processes represented as stochastic It\^o integrals (possibly, with a drift term). Usually, these integrals have unlimited range of possible values.
However, there are special cases of integrals with limited range. These integrals can be obtained, for instance, as conditional expectations of random variables with
limited  range, or via restriction of the integration interval by a random Markov stopping times  preventing the range growth. These approaches may be inconvenient in some cases.
For example, calculation of a condition expectation is essentially a forecast of a contingent claim depending on the future values, and this procedure can be difficult.
Besides, one would need to specify first this contingent claims.  On the other hand,   restriction of the integration interval by stopping times
leads to stochastic integrals with some paths being frozen at that stopping times. Obviously, this feature  could be undesirable.

The present paper  suggests special examples of stochastic It\^o integrals with controlled growth  of their containing range.  The integrands for this integrals are presented
explicitly. The paper uses an original  approach does not involve neither stopping times nor forecasting or calculation of the conditional expectations of a contingent claim. This approach does not involve neither forecasting nor calculation of the conditional expectations of a contingent claim.

\section{The main result}\label{SecMain}
We are also given a standard  complete probability space $(\O,\F,\P)$ and a
right-continuous filtration $\{\F_t\}_{t\ge 0}$ of complete $\s$-algebras of
events. In addition, we are given  an one-dimensional Wiener
process $w(t)|_{t\ge 0}$, that is a Wiener
process adapted to  $\{\F_t\}$ and such that $w(0)=0$ and that $\F_t$ is independent from $w(s)-w(q)$ if $t\ge s>q\ge 0$.

Consider a  continuous time one-dimensional  random
process $x(t)|_{t\ge 0}$
such that
\baaa
dx(t)=a(t)dt+\s(t) dw(t).
\eaaa
Here  $a(t)$  and $\s(t)$  are bounded real-valued one-dimensional $\F_t$-adapted processes.

Let  $u(t):(0,+\infty)\to \R$ be a real valued random $\F_t$-adapted process
that is  integrable on any finite time interval. For simplicity, we assume that $\esssup_{\o\in\O}\int_0^t|u(s)|ds<+\infty$ for any $t>0$.

Let $\{\GG_t\}_{t\ge 0}$ be the  filtration of complete $\s$-algebras of
events generated by the process $(x(t),u(t))$.

It can be noted that, since  the process $\s(t)$ is adapted
to the filtration generated by $x(t)$,  it follows that $\{\GG_t\}_{t\ge 0}$ is also the  filtration  generated by the process $(x(t),u(t),\s(t))$; see e.g. Remark 1.1 in \cite{D02}, p.10, or
 Proposition 7.1 in \cite{D07}, where this was shown for a log-normal types
 of processes which was rather technical.
\begin{theorem}\label{Th1}
Consider  processes $X(t)$ and $Y(t)$ defined for $t\in[0,+\infty)$ as
\baa
X(t)=\int_0^t \cos(x(t)-x(s))u(s)ds
\label{X}
\eaa
 and \baa &&Y(t)= \int_0^t \sin(y(t)-y(s))u(s)ds.\label{Yint}\eaa
In this case,  \baa
&&\int_0^t X(s) dx(s) =Y(t)+\frac{1}{2}\int_0^t \s(s)^2Y(t) ds.
\label{Y}\eaa
Then  \baaa
 \sqrt{X(t)^2+Y(t)^2}\le \int_0^t|u(s)|dt\quad \hbox{a.s.}\quad \forall a(\cdot).
\eaaa
\end{theorem}
\par
Clearly, the process  $X(t)$  is bounded uniformly in all $a(\cdot)$ almost surely on any finite time interval and $\GG_t$-adapted. Hence the stochastic integral (\ref{X})
is well defined and is  bounded uniformly in all $a(\cdot)$ almost surely on any finite time interval.

Representation (\ref{Yint}) in Theorem \ref{Th1} implies that the boundaries for the range of the stochastic integral $Y(t)$ are defined by the choice of the process $u$.
Respectively, Theorem \ref{Th1} allows to construct stochastic processes with preselected  on a given time interval
time depending boundaries for their  range.

\subsection{Proof  of Theorem \ref{Th1}}
The proof follows the idea of the proof of
Lemma 3.2 from \cite{D14} (see also \cite{D18}).

Consider a process
 \baa &&dZ(t)=iZ(t) dx(t)-\frac{1}{2}\s(t)^2Z(t)dt + u(t)dt,\quad
t\in(0,\infty),\nonumber\\&& Z(0)=0.\label{Z}\eaa

In (\ref{Z}), $i=\sqrt{-1}$ is the imaginary unit.

\begin{lemma}\label{lemmaZ}
For any $T>0$, we have that
\baaa
Z(t)=i\int_0^t e^{i [x(t)-x(s)]}u(s)ds
\eaaa
and
\baaa
\Im Z(t)=\int_0^t \Re Z(s) dx(s)\quad \forall t>0.
\eaaa
The process $Z(s)$  is  $\GG_t$-adapted and such that
\baa
 |Z(t)|\le \int_0^t |u(s)|ds\quad  \forall t>0.
\label{Zest}\eaa
\end{lemma}

{\em Proof of Lemma \ref{lemmaZ}.}
Let $F(t,s)$ be defined as the solution of the It\^o equation
 \baaa  &&dF(t,s)=iF(t,s)[a(t)dt+\s(t)
dw(t)]-\frac{1}{2}F(t,s)\s(t)^2
dt, \quad t>s\ge 0,\qquad\\&& F(s,s)=1.\label{CT22} \eaaa
By the It\^o formula,
\baaa F(t,s)&=&F(s,s)\exp\left(i\int_{s}^{t} a(r)dr
+i\int^{t}_{s}\s(r) dw(r)-\frac{i^2}{2}\int^{t}_{s}\s(r)^2dr-\frac{1}{2}\int^{t}_{s}\s(r)^2dr \right)\nonumber\\&=&\exp\left(i\int_{s}^{t}a(r)dr+i\int^{t}_{0}\s(r) dw(r)
\right)\nonumber\\
&=&\exp\left(i[y(t)-y(s)]\right)
\quad\hbox{a.s.}. \label{solSx} \eaaa In particular, we have that
\baaa |F(t,s)|=1 \quad\hbox{a.s.}\label{SYxx} \eaaa

Direct differentiation gives that
\baaa
Z(t)=\int_0^t F(t,s)  u(s)ds.
\eaaa
Hence
 \baaa
|Z(t)|\le \int_0^t |F(t,s)| |u(s)|ds\le \int_0^t  |u(s)|ds .
\label{Zest2}\eaaa
Let $X(t)=\Re Z(t)$ and  $Y(t)=\Im Z(t)$.Then (\ref{X}) and (\ref{Yint}) hold.
 Further, we have from (\ref{Z}) that
\baaa
&&dX(t)=-Y(t)dx(t)+u(t)dt-\frac{1}{2}\s(t)^2 X(t)dt,\\
&&dY(t)=X(t)dx(t)-\frac{1}{2}\s(t)^2 Y(t)dt.
\eaaa
Then the proof of Lemma \ref{lemmaZ}  follows. $\Box$

{\em Proof of Theorem \ref{Th1} } follows from this and
from (\ref{Zest}). $\Box$
 \begin{remark}\label{rem1}{\rm
 Consider the case where $a(t)\equiv 0$.
 In this case, it follows from the proof above that, for any $\F_t$-adapted
 process $\s(\cdot)$,  there exists a $\F_t$-adapted process $U:[0,T]\times\O \to \C$ such that
 $|U(t)|=1$ and that
 the integral
$\int_0^T\s(t) U(t)dw(t)$ has a limited range in $\C$. To see this, it suffices to
select $U(t)=F(t,0)$ and observe that $i\int_0^T\s(t) F(t,0)dw(t)=F(T,0)-1 +\frac{1}{2}
\int_0^T\s(t)^2 F(t,0)dt$ and that
\baaa
\left|F(T,0)-1 +\frac{1}{2}
\int_0^T\s(t)^2 F(t,0)dt\right|\le
\left|F(T,0)\right|+1 +\frac{1}{2}
\int_0^T\s(t)^2 |F(t,0)|dt
\\ \le 2 +\frac{1}{2}
\int_0^T\s(t)^2 dt.
\eaaa
 }
 \end{remark}
\section{Some modifications}
The approach demonstrated above allows many modifications.
Let us provide one of possible  modifications,

\begin{theorem}\label{DTh1}
Consider a process $\ww Y(t)$ defined as
\baa
\ww Y(t)=\int_0^t \cos(x(t)-x(s))u(s)e^{\frac{1}{2}\int_s^t\s(r)^2dr}ds,
\label{DY}
\eaa
where $t\in[0,+\infty)$. Further, let a process $\ww X(t)$ be defined as the stochastic integral  \baa
&&\ww X(t)=-\int_0^t \ww Y(s) dx(s),
\label{DX}\eaa
  Then $\ww X(t)$ can be represented as
 \baa &&\ww X(t)= -\int_0^t \sin(y(t)-y(s))e^{\frac{1}{2}\int_s^t\s(r)^2dr}u(s)ds.\label{DXint}\eaa
\end{theorem}
\par
Clearly, the process  $\ww Y(t)$  is bounded uniformly in all $a(\cdot)$ almost surely on any finite time interval and $\GG_t$-adapted. Hence the stochastic integral (\ref{DX})
is well defined.

Representation (\ref{DXint}) in Theorem \ref{DTh1} implies that the boundaries for the range of the stochastic integral $\ww X(t)$ are defined by the choice of the process $u$.
Respectively, Theorem \ref{DTh1} allows to construct stochastic processes with preselected  on a given time interval
time depending boundaries for their  range.
\begin{corollary}\label{DThC} Let $T>0$ be fixed, and let $\psi(t):(0,T)\to \R$ be a  integrable function.
Let $u(t)=e^{-\frac{1}{2}\int_t^T\s(s)^2ds}\psi(t)$.
Then \baa
&& \ww X(t)= -\int_0^t \sin(y(t)-y(s))\psi(s)ds,
\label{DX2}\\
&&\ww Y(t)=\int_0^t \cos(y(t)-y(s))\psi(s)ds,\quad t\in[0,T],\quad
T\in(0,\infty),\label{DY2}\eaa
and  \baaa
|\ww X(t)|\le \sqrt{\ww X(t)^2+\ww Y(t)^2}\le \int_0^t|\psi(s)|dt\quad \hbox{a.s.}\quad \forall a(\cdot).
\eaaa
\end{corollary}
\begin{example}
\begin{enumerate}
\item
if $|\psi(t)|=\a^{-1}t^{\a-1}$ for $\a >1/2$,
then $\sqrt{\ww X(t)^2+\ww Y(t)^2}\le t^{\a}$ a.s. for all $a(\cdot)$.
\item If $|\psi(t)|=1/(a+t)$ for some $q>0$,
then $ \sqrt{\ww X(t)^2+\ww Y(t)^2}\le \ln(q+t)-\ln q$ a.s. for all $a(\cdot)$.
\end{enumerate}
\end{example}

\subsection{Proof  of Theorem \ref{DTh1}}
The proof is similar to the proofs of Theorem \ref{Th1}; we provide it for completeness.

Consider a process
 \baa &&d\ww Z(t)=i[\ww Z(t) dx(t) + u(t)dt],\quad
t\in(0,\infty),\nonumber\\&& \ww Z(0)=0.\label{DZ}\eaa

In (\ref{DZ}), $i=\sqrt{-1}$ is the imaginary unit.

\begin{lemma}\label{DTh11}
For any $T>0$, we have that
\baaa
\ww Z(t)=i\int_0^t e^{i [x(t)-x(s)]}e^{\frac{1}{2}\int_s^t\s(r)^2dr}u(s)ds
\eaaa
and
\baaa
\Re \ww Z(t)=\int_0^t \Im \ww Z(s) dx(s)\quad \forall t>0.
\eaaa
The process $\ww Z(s)$  is  $\GG_t$-adapted and such that
\baaa
 |\ww Z(t)|\le e^{\frac{1}{2}\int_0^T\s(s)^2ds}  |u(t)|\quad  \forall T>0, \quad\forall t\in[0,T].
\eaaa
\end{lemma}

{\em Proof of Lemma \ref{DTh1}.}
Let $\ww F(t,s)$ be defined as the solution of the It\^o equation
 \baaa  d\ww F(t,s)=i\ww F(t,s)[a(t)dt+\s(t)
dw(t)], \quad t>s\ge 0,\qquad \ww F(s,s)=1.\label{DCT22} \eaaa
By the It\^o formula,
\baaa \ww F(t,s)&=&\ww F(s,s)\exp\left(i\int_{s}^{t} a(r)dr
+i\int^{t}_{s}\s(r) dw(r)-\frac{i^2}{2}\int^{t}_{s}\s(r)^2d\right)\nonumber\\&=&\exp\left(i\int_{s}^{t}a(r)dr+i\int^{t}_{0}\s(r) dw(r)+\frac{1}{2}\int^{t}_{s}\s(r)^2dr
\right)\nonumber\\
&=&\exp\left(i[y(t)-y(s)]+\frac{1}{2}\int^{t}_{s}\s(r)^2dr\right)
\quad\hbox{a.s.}. \label{DsolSx} \eaaa In particular, we have that
\baaa |\ww F(t,s)|=\exp\left(\frac{1}{2}\int^{t}_{s}\s(r)^2dr\right) \quad\hbox{a.s.}\label{DSYxx} \eaaa

Direct differentiation gives that
\baaa
\ww Z(t)=i\int_0^t \ww F(t,s)  u(s)ds.
\eaaa
Hence
 \baaa
|\ww Z(t)|\le \int_0^t |\ww F(t,s)| |u(s)|ds\le \int_0^t e^{\frac{1}{2}\int^t_t\s(r)^2 dr} |u(s)|ds .
\eaaa

Let $\ww X(t)=\Re \ww Z(t)$ and  $\ww Y(t)=\Im \ww Z(t)$. We have that
\baaa
&&d\ww X(t)=-\ww Y(t)dx(t),\\
&&d\ww Y(t)=\ww X(t)dx(t)+u(t).
\eaaa
It follows that
\baaa
\ww XT)=\int_0^T\g(t)dx(t),
\eaaa
where   $\g(t)=-\ww Y(t)$.  We have that \baaa
\int_0^t|\g(s)|ds\le \int_0^t|\ww Z(s)|ds \le \int_0^s e^{\frac{1}{2}\int_0^T\s(r)^2 dr}  |u(s)|ds
\le  e^{\frac{1}{2}\int_0^T\s(r)^2 dr}  \int_0^t |u(s)|ds\quad \forall T>0,\quad t\in[0,T].
\eaaa
Hence $|\g(t)|\le |\ww Z(t)|\le e^{\frac{1}{2}\int_0^T\s(r)^2dt}|u(t)|$ for any $T>0$ and $t\in[0,T]$.

Then the proof of Lemma follows. $\Box$

The proof of Theorem \ref{DTh1} follows from the lemma. The proof of  Corollary
 \ref{DThC}  follows from the theorem applied to the corresponding choice of $u$.

 \begin{remark}{\rm  For the case where $a(t)\equiv 0$,
similarly Remark \ref{rem1}, it can be shown that, that, for any $\F_t$-adapted
 process $\s(\cdot)$,  there exists a $\F_t$-adapted process $U:[0,T]\times\O \to \C$ such that
 $|U(t)|=\exp\left(\frac{1}{2}\int^{t}_{o}\s(r)^2dr\right)$ and that
 the integral
$\int_0^T\s(t) U(t)dw(t)$ has a limited range in $\C$. To see this, it suffices to
select $U(t)=i\ww F(t,0)$ and observe that $\int_0^T\s(t) i \ww F(t,0)dw(t)=\ww F(T,0)-1$. }
 \end{remark}
\section{Possible  applications for financial modelling}
One of core problem of financial mathematics is the portfolio selection problem.
Application of classical methods of optimal stochastic control for portfolio optimization problems
  requires
 forecasting of market parameters.
This forecasting is usually  difficult.
This problem is related to the open  problem of
validation of the so-called {\em technical analysis} methods that offer trading
strategies based on historical observations.  There are many
different strategies suggested in this  framework (see, e.g., \cite{B,BG,C,Hsu1,Hsu,D02,DS4}
and the references therein. It is
 known that  mean-reverting market models and market models with bounded range
 for the prices generate
 some special speculative
opportunities (see, e.g., \cite{CF,C,D06,D7,D12,ff,Lo,Lor,Shi}).
Theorems \ref{Th1}-\ref{DTh1}  give a possibility to convert a stock price process $x(t)$
into a processes $Y(t)$ or  $X(t)$ that could have features similar to mean-reverting market models and market models with bounded range for the prices.
For this new artificial asset, one can apply strategies \index{defined in Theorem 2.1} from \cite{CF,C,D06,D7,D12,ff}.
We leave this for the future research.

\begin{figure}[ht]
\centerline{\psfig{figure=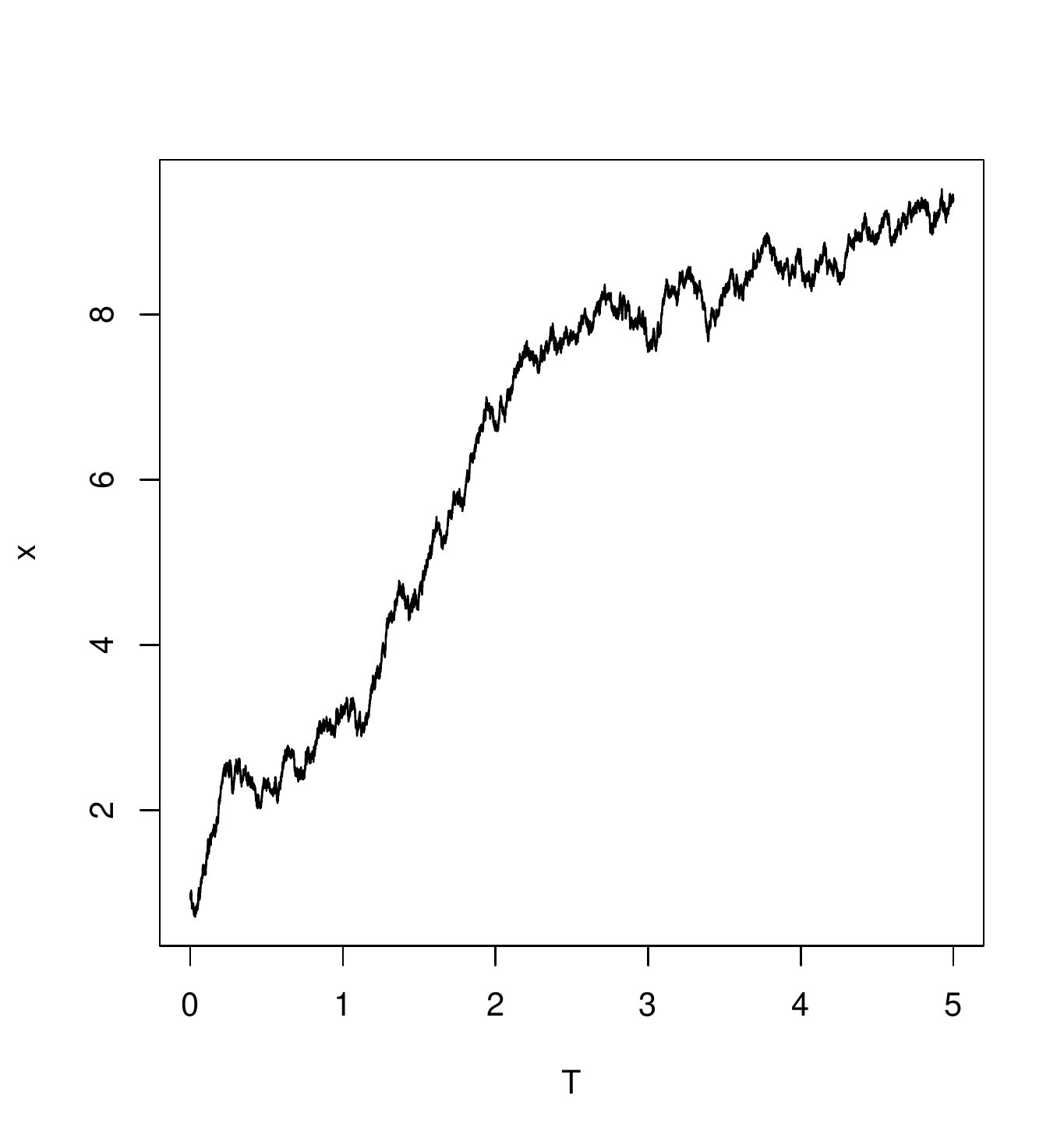,width=9cm,height=6.0cm}}
\caption[]{\sm
 $x(t)$.
 }\label{figx}
  \centerline{\psfig{figure=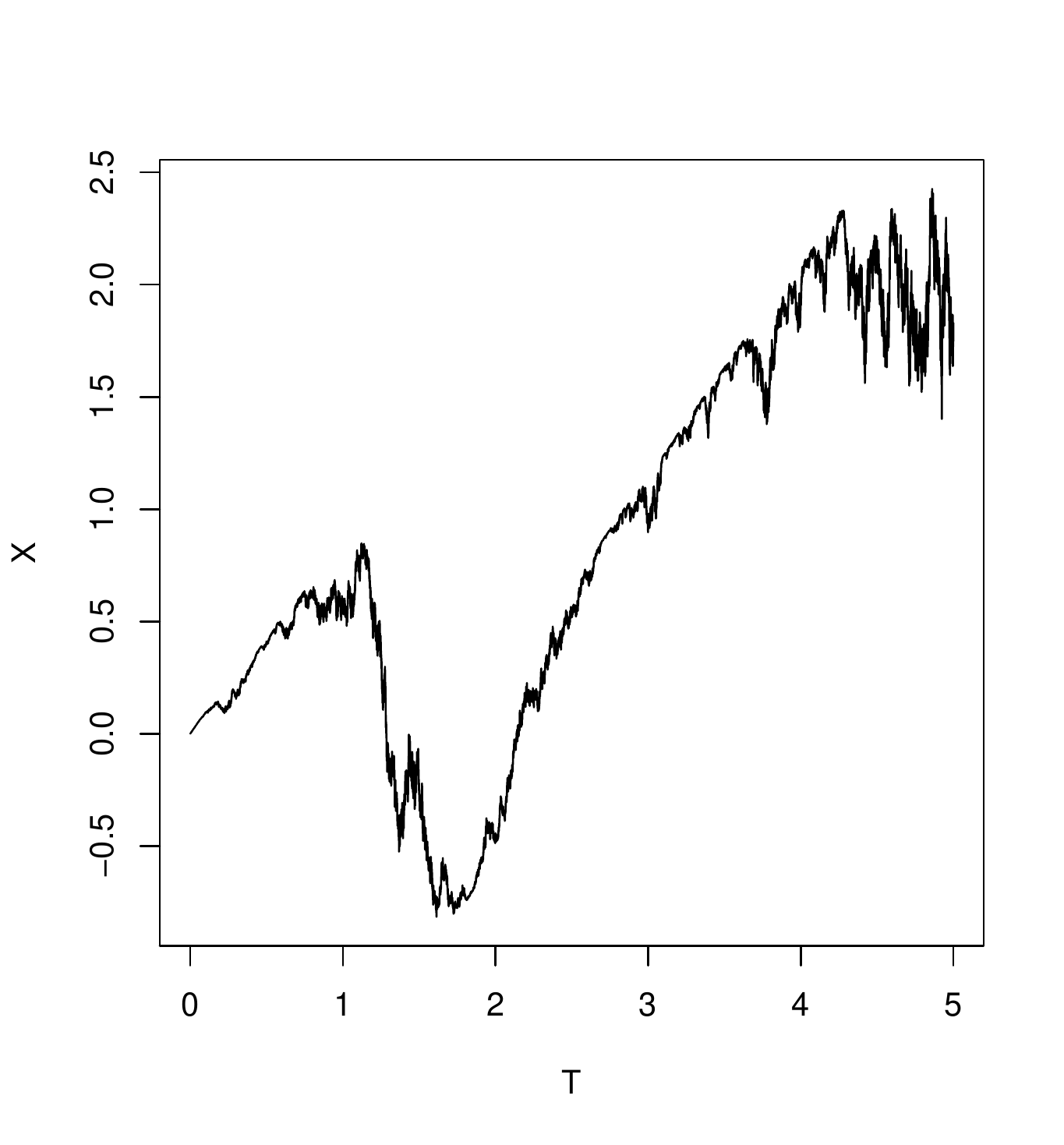,width=9cm,height=6.0cm}}
\caption[]{\sm
 $X(t)$.}\label{figX}
 \centerline{\psfig{figure=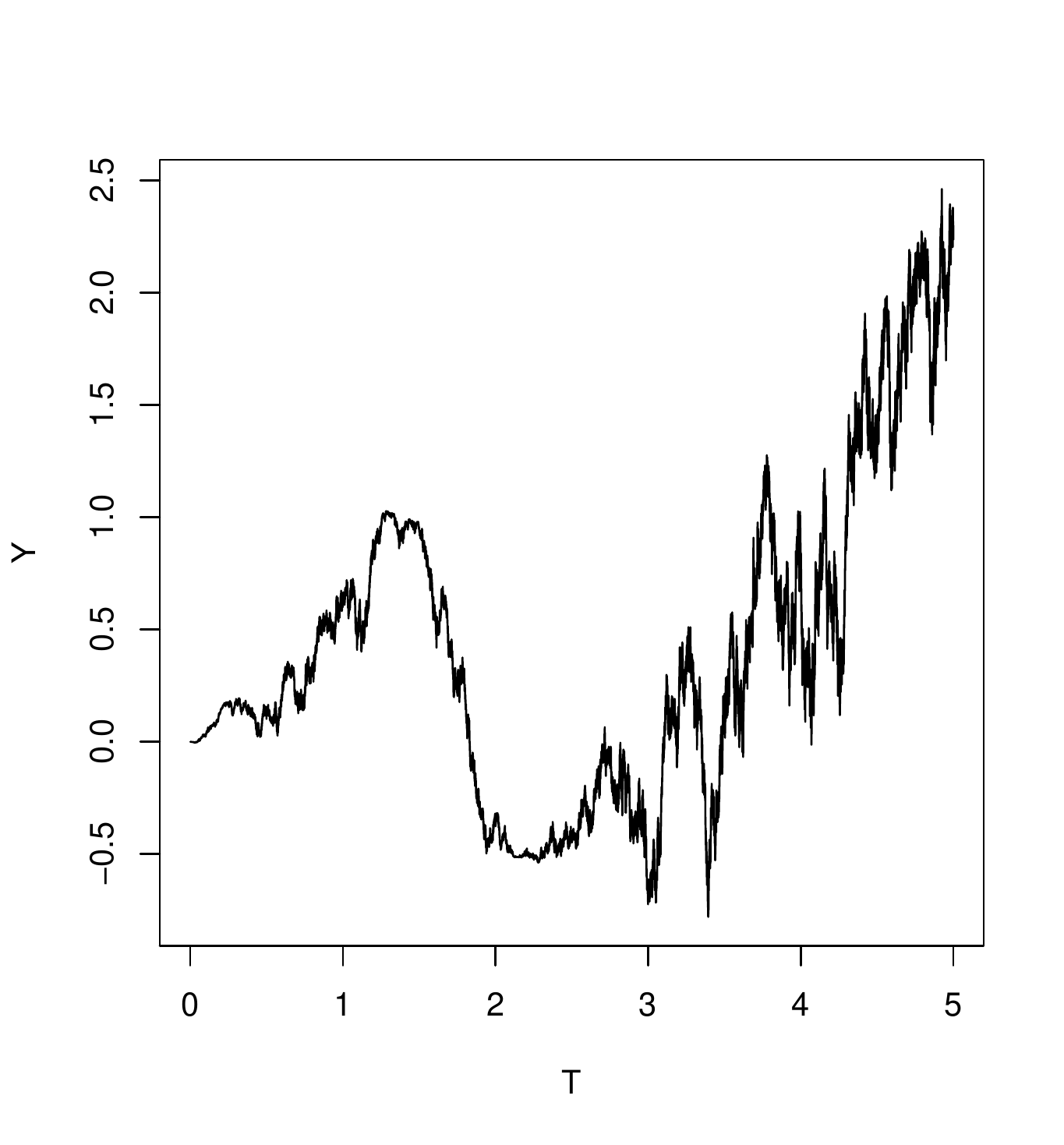,width=9cm,height=6.0cm}}
\caption[]{\sm
 $Y(t)$.}\label{figY}
  \end{figure}

\begin{figure}[ht]
 \centerline{\psfig{figure=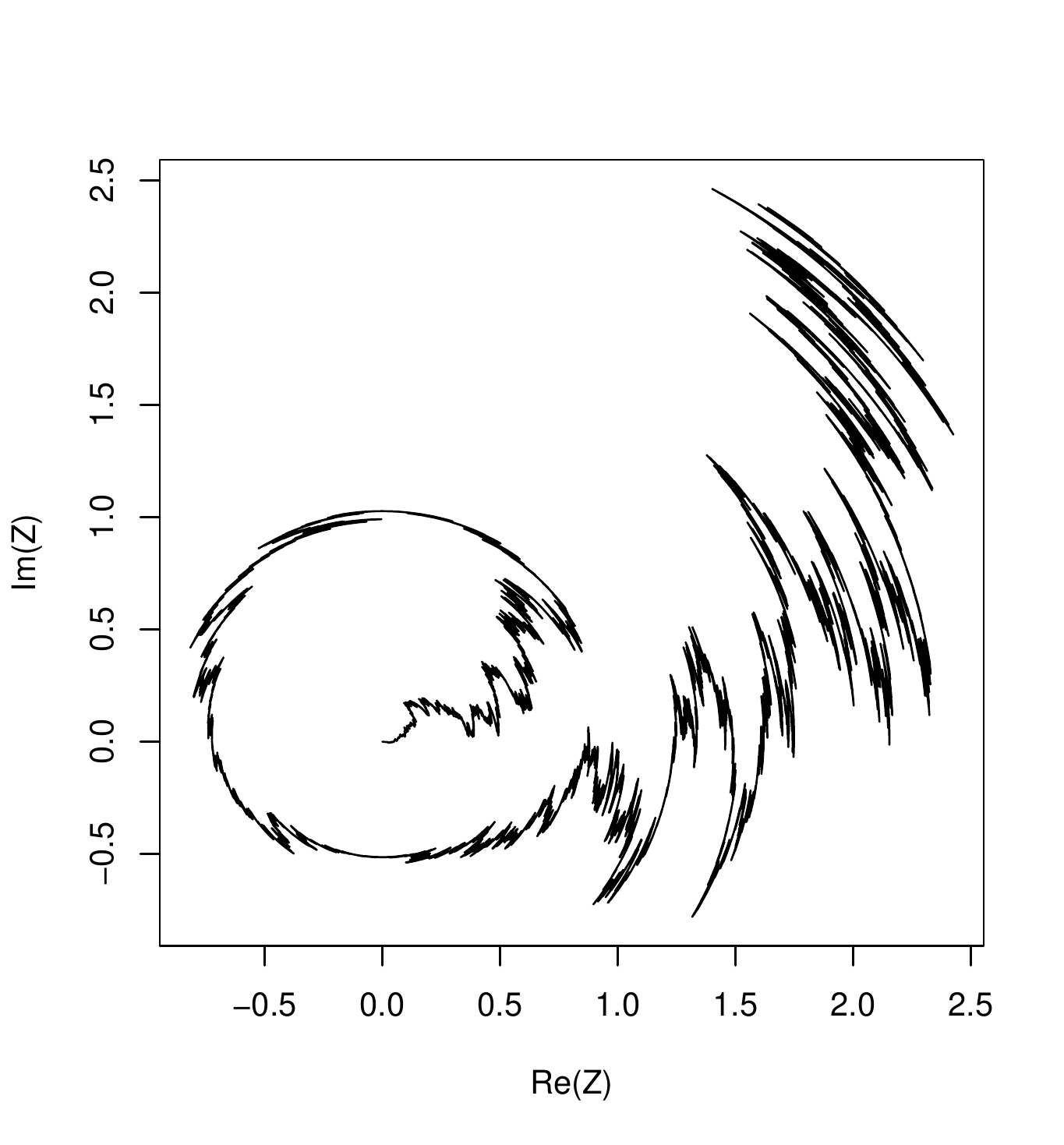,width=6cm,height=6.0cm}}
\caption[]{\sm
 $Z(t)=X(t)+i Y(t)$.}
 \centerline{\psfig{figure=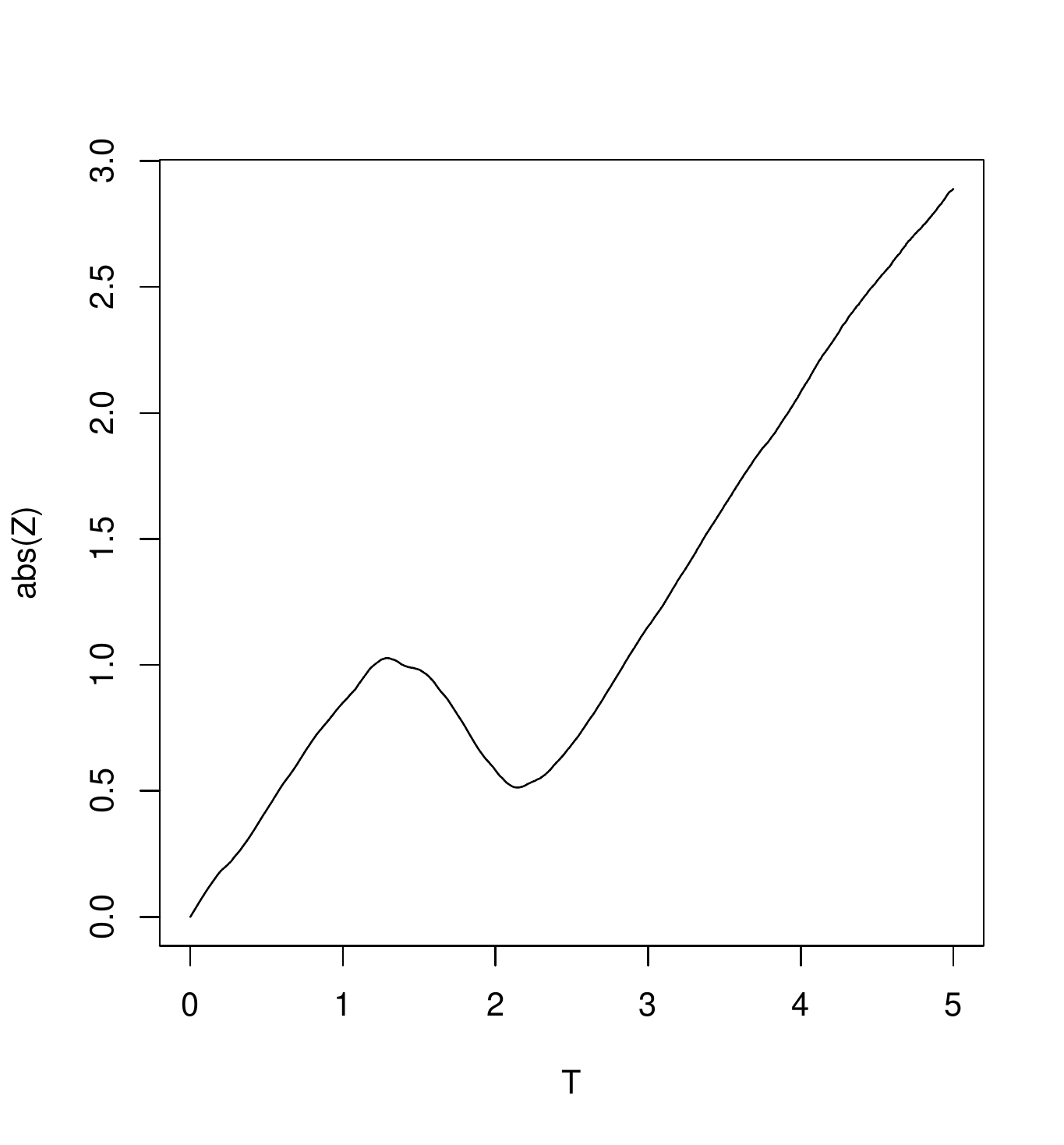,width=9cm,height=6.0cm}}
\caption[]{\sm
 $|Z(t)|$.}
 \vspace{0cm}
 \centerline{\psfig{figure=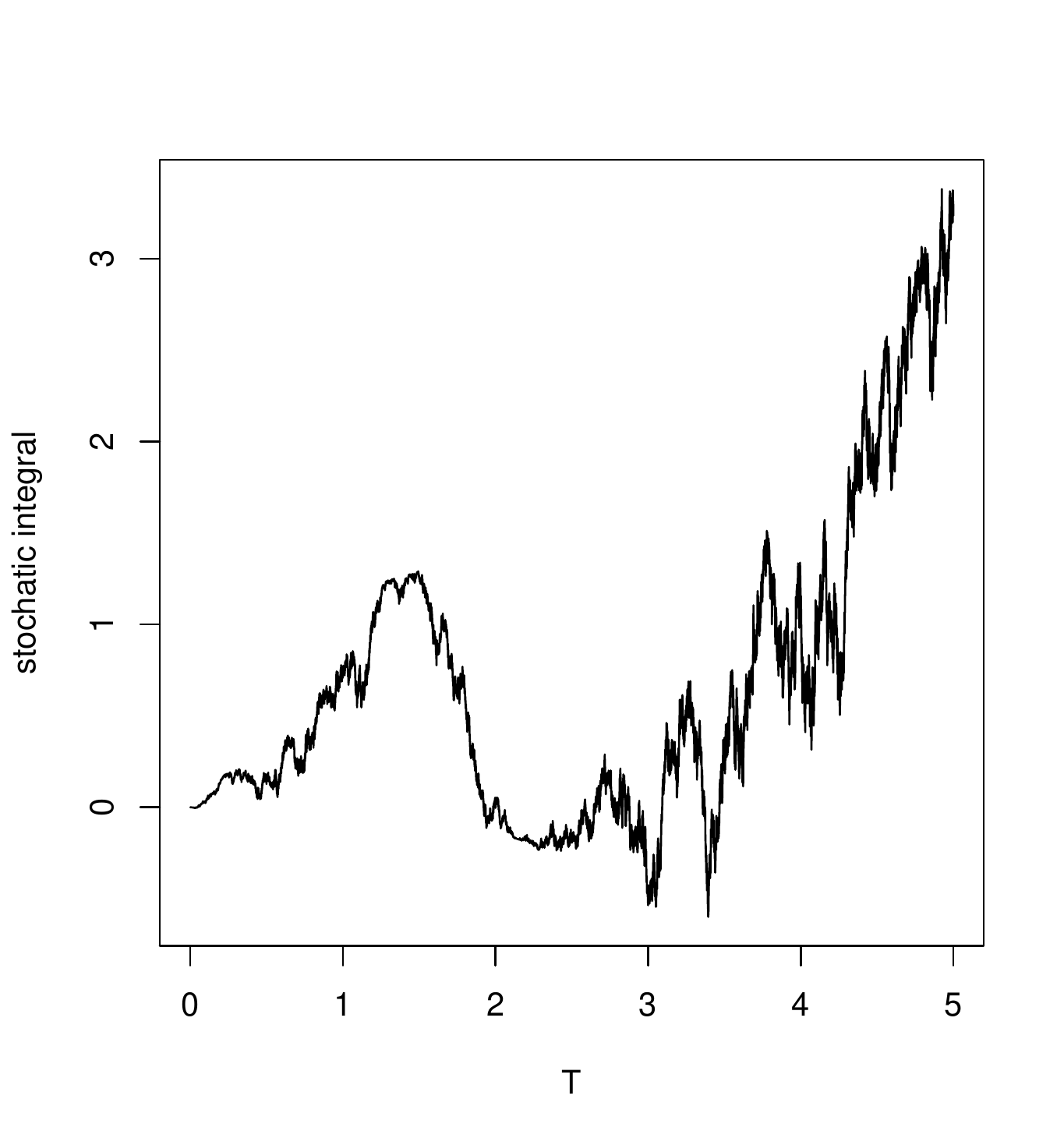,width=9cm,height=6.0cm}}
\caption[]{\sm
 $\int_0^tX(s)dx(s)$.}\label{figF}
 \vspace{0cm}\end{figure}

Figures \ref{figx}-\ref{figF} shows sample paths of processes
introduced above and obtained via Monte-Carlo simulation
under the assumption that $T=5$, $\s(t)\equiv\s= 1$, $a(t)\equiv a=2$, and $u(t)\equiv 1$.
We used natural discretization in time with $10^5$ grid point on the interval $[0,T]$.
Calculations were executed using R and RStudio programms.

\end{document}